\documentclass[amsfonts,amssymb]{article} 
\usepackage{amssymb}
\usepackage{graphics}
\usepackage{epsfig}
\title{Fundamental Solutions for  Hyperbolic  Operators 
  with Variable Coefficients}

\author{\textsc{Karen Yagdjian}}

\date{} 
 
\newtheorem{theorem}{Theorem}[section]
\newtheorem{corollary}{Corollary}

\newtheorem{definition}[theorem]{Definition}

\begin{document}
\maketitle

\thispagestyle{empty} 
\begin{center}
Department of Mathematics, \\
University of Texas-Pan American,  \\
1201 W.~University Drive, \\
Edinburg, TX 78541-2999,  
USA \\
yagdjian{@}utpa.edu 
\end{center}

\begin{abstract}
{In this   article  we describe the novel method to construct fundamental solutions for operators with variable coefficients.
That method was  introduced   in \cite{YagTricomi} to study the Tricomi-type equation. 
More precisely, the new integral  operator  is suggested which transforms the family of the fundamental solutions of   
 the Cauchy problem for the equation with the constant coefficients  to the fundamental solutions for the 
operators with variable coefficients. 
 }\end{abstract}

\section{Introduction}
\label{S1}

\setcounter{equation}{0}
\renewcommand{\theequation}{\thesection.\arabic{equation}}

In this  article  we describe the novel method to construct fundamental solutions for operators with variable coefficients.
We also give a brief survey  of some results obtained by that method, which   was  introduced   in \cite{YagTricomi} to solve the Cauchy problem for the Tricomi-type equation. 
Later on it was applied to several partial differential equations 
with variable coefficients containing  
some equations arising in the mathematical cosmology.   
More precisely, the new integral operator is suggested which transforms the family of the fundamental solutions of   
 the Cauchy problem for the equation with the constant coefficients  to the fundamental solutions for the 
operators with variable coefficients.  The kernel of that transformation  contains Gauss's hypergeometric function. 
\smallskip

This method  was  used in \cite{Galstian_Kinoshita_Yagdjian,Kinoshita_Yagdjian},\cite{YagTricomi}-\cite{
yagdjian_DCDS}   
to investigate in the unified way several equations such as 
the linear  and semilinear Tricomi and Tricomi-type equations, Gellerstedt equation, the wave equation in Einstein-de~Sitter spacetime, the wave and the Klein-Gordon 
equations in the de~Sitter and anti-de~Sitter spacetimes.
The listed  equations play important role in the gas dynamics, elementary particle physics, quantum field theory in the curved spaces, and cosmology.
 For all above mentioned equations, we have obtained among  other things, fundamental solutions, the representation formulas for the initial-value problem, the $L_p-L_q$-estimates, local and global solutions for the semilinear equations, blow up phenomena,
self-similar solutions and number of other results.  
\smallskip

The starting point of our approach is the Duhamel's principle, which we revise in order to prepare the ground for generalizations. 
It is well-known that the solution of the Cauchy problem for the string equation with the source term $f$,  $f(x,t) \in C^\infty({\Bbb R}^{2})$, 
\begin{equation}
\label{main_n=1}  
 u_{tt}-u_{xx}=f(x,t) \quad \mbox{\rm in}  \,\, {\mathbb R}^{2}, \quad
u(x,t_0)=0,\quad  u_t(x,t_0)=0   \quad \mbox{\rm in}  \,\, {\mathbb R} \,, 
\end{equation}
can be written as an integral 
\[ 
u(x,t)= \int_{t_0}^t v(x,t;\tau ) \,d\tau 
\]
of the family of the solutions $ v(x,t;\tau ) $ of the problem without the source term, but with the second initial datum 
\[ 
 v_{tt}-v_{xx}= 0\quad \mbox{\rm in}  \,\, {\mathbb R}^{2}, \quad
v(x,\tau ;\tau )=0,\quad  v_t(x,\tau ;\tau )=f(x,\tau )   \quad \mbox{\rm in}  \,\, {\mathbb R}\,. 
\]

Our {\it first observation} is that we obtain the following representation of the solution of (\ref{main_n=1})
\begin{equation}
\label{main}  
u(x,t)= \int_{t_0}^t \,d\tau \int_{ 0}^{  t-\tau  } w(x,z;\tau )\,dz\,,  
\end{equation}
if we denote
\[
w(x;t;\tau ):=  \frac{1}{2}\left[ f(x+t,\tau )+ f(x-t,\tau )\right],
\]
where the function $w=w(x;t;\tau ) $ is the solution of the problem
\begin{eqnarray} 
\label{02new}
w_{tt}-w_{xx}= 0\quad \mbox{\rm in}  \,\, {\mathbb R}^{2}, \quad 
w(x,0 ;\tau )=f(x,\tau ),\quad  w_t(x,0 ;\tau )= 0  \quad \mbox{\rm in}  \,\, {\mathbb R} \,.  
\end{eqnarray}
This formula allows us to solve problem with the source term if we solve the problem for the {\sl same} equation without source 
term but with the first initial datum.
We claim that the  formula (\ref{main}) can be used also for the wave equation with $x \in {\mathbb R}^{n}$, for all $n \in {\mathbb N}$. (See, e.g, \cite{YagTricomi}.)
More precisely, it holds also for the problem
\[ 
 u_{tt}-\Delta u =f(x,t) \quad \mbox{\rm in}  \,\, {\mathbb R}^{n+1}, \quad 
u(x,t_0)=0,\quad  u_t(x,t_0)=0   \quad \mbox{\rm in}  \,\, {\mathbb R}^{n}\,,  
\]
with the function $w=w(x;t;\tau ) $ solving 
\begin{equation}
\label{wn}   
  w_{tt}-\Delta w = 0\,\, \mbox{\rm in}  \,\, {\mathbb R}^{n+1}, \quad
w(x,0 ;\tau )=f(x,\tau ),\,\, w_t(x,0 ;\tau )= 0  \,\, \mbox{\rm in}  \,\, {\mathbb R}^{n}.  
\end{equation}
Note that in the last problem the initial time $t=0$ is frozen, while in the Duhamel's principle it is varying with the parameter $ \tau $. 
 
The  {\it second observation} is that in (\ref{main}) the upper limit $t-\tau $  of the inner integral is generated by the  propagation phenomena  with the speed which is equals to one. 
In fact, that is a distance function between the points at time $t$ and $\tau $. 

Our {\it third observation}  is that the solution operator $G\,:\,f \longmapsto u $ can be regarded as a composition of two operators. The first one 
\[
{\mathcal WE}: \,\, f  \longmapsto  w 
\]
is a Fourier Integral Operator (FIO), which is a solution operator of the Cauchy problem with the first initial datum for wave equation in the Minkowski spacetime.  The second operator  
\[
{\mathcal K}:\,\,w\longmapsto   u 
\]
is the integral operator given by (\ref{main}). We regard the variable $z$ in  (\ref{main}) as a ``subsidiary time''. Thus, $G= {\mathcal K}\circ {\mathcal WE}$ and we arrive at the diagram: 
%\begin{eqnarray*}
%&    {\mathcal  WE}  & \\
%\vspace{-1cm}  f &  \longmapsto     & w  \\
%  &\searrow G \hspace{0.7cm}   {\mathcal K} \swarrow&  \\
% & u &
% \end{eqnarray*}
\vspace{-0.3cm}
\begin{figure}[h]
\begin{center}
\scalebox{0.15}{\includegraphics{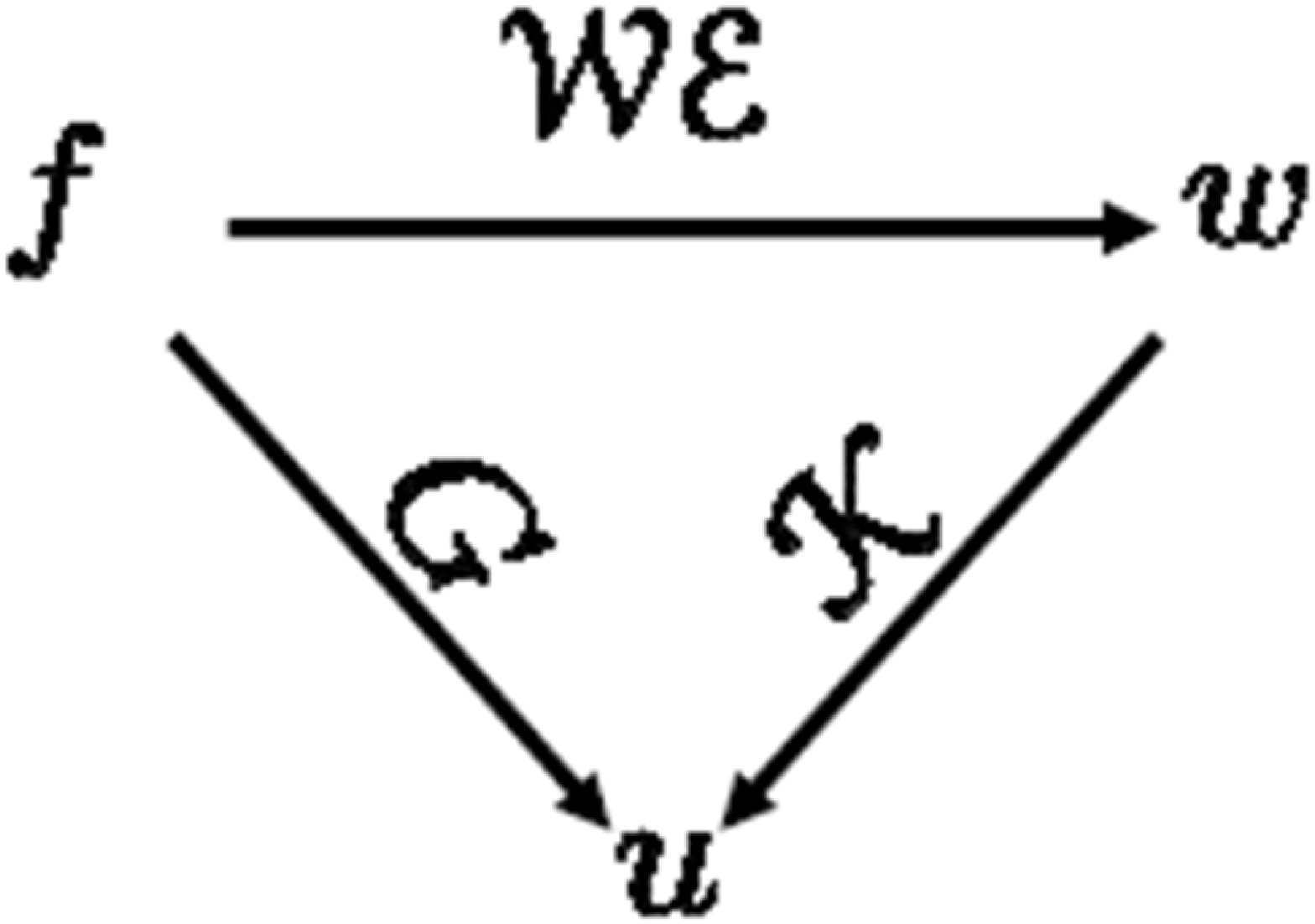}}
\label{F1}
\end{center}
\end{figure}
\vspace{-0.5cm}

Our aim is that based on this diagram to generate a class of operators for which we will obtain   explicit representation formulas for the solutions. 
That means also that we will have representations for the fundamental solutions of the partial differential operator. In fact,  this diagram
brings into a single  hierarchy several different partial differential operators.
Indeed, if we take into account the propagation cone by introducing the distance function 
$ \phi (t)$, and if we provide   the integral operator with the kernel $ K (t;r,b) $ as follows:
\begin{equation} 
\label{Oper_K}
{\mathcal K} [w](x,t)  
 = 
2   \int_{ t_0}^{t} db 
  \int_{ 0}^{ |\phi (t)- \phi (b)|}   K (t;r,b)  w(x,r;b )  dr  
, \quad x \in {\mathbb R}^n, \,\, t>t_0, 
\end{equation}
then we actually can generate new representations for the solutions of the different well-known equations. Below we illustrate the suggested scheme 
by several examples. 
\smallskip

\noindent
$1^0 $ Klein-Gordon equation in  the Minkowski spacetime. 
If we choose the kernel $K (t;r,b)$ as
\begin{eqnarray} 
\label{0.6}
K (t;r,b)  
& = &
J_0\left(\sqrt{(t-b)^2-r^2}\right)\,, 
\end{eqnarray}
where $J_0(z) $ is the  Bessel function of the first kind, and if we choose the distance function as $\phi (t)=t $, 
then we can  prove (see Theorem~\ref{T0_1} below) that the function  
\[
u(x,t)  
 =  
\int_{ t_0}^{t} db 
  \int_{ 0}^{t-b }   J_0\left(\sqrt{(t-b)^2-r^2}\right)   w(x,r;b )  dr  
, \,\, x \in {\mathbb R} , \,\, t>t_0, 
\]
solves the problem for the Klein-Gordon equation with a positive mass equals to $ 1$ in the one-dimensional Minkowski spacetime,
\begin{eqnarray*} 
 u_{tt}-u_{xx}+u =f(x,t) \quad \mbox{\rm in}  \,\, {\mathbb R}^{2}, \quad
u(x,t_0)=0,\quad  u_t(x,t_0)=0   \quad \mbox{\rm in}  \,\, {\mathbb R} \,,  
\end{eqnarray*}
provided that $ w(x,r;b )$ is a corresponding solution of the problem for the wave equation in the Minkowski spacetime.
We emphasis that the function  $w= w(x,t;b )$, with  $b$ regarded as a parameter, and  the function $ u = u(x,t)$ solve {\sl different} equations.
This is a fundamental distinction   from the Duhamel's principle.

Now if we  choose the kernel $K (t;r,b)$ as
\begin{eqnarray} 
\label{BesselI_0}
K (t;r,b)  
& = &
I_0\left(\sqrt{(t-b)^2-r^2}\right) \,,
\end{eqnarray}
where $I_0(z) $ is the modified Bessel function of the first kind, and the distance function as $\phi (t)=t $, 
then the function  
\[ 
u(x,t)  
=  
\int_{ t_0}^{t} db 
  \int_{ 0}^{t-b }   I_0\left(\sqrt{(t-b)^2-r^2}\right)  w(x,r;b )  dr  
, \,\, x \in {\mathbb R}^n, \,\, t>t_0, 
\]
solves the problem for the Klein-Gordon equation  with an imaginary mass  in  the  one-dimensional Minkowski spacetime,
\begin{eqnarray*} 
u_{tt}-u_{xx}-u =f(x,t) \quad \mbox{\rm in}  \,\, {\mathbb R}^{2},  \quad
u(x,t_0)=0,\quad  u_t(x,t_0)=0   \quad \mbox{\rm in}  \,\, {\mathbb R} \,, 
\end{eqnarray*}
provided that $ w(x,r;b )$ is a corresponding solution of the problem (\ref{02new}) for the wave equation in the one-dimensional Minkowski spacetime. 
According to the next theorem the representation formulas valid also for the higher dimensional equations. The proof   is by straightforward substitution. 
\begin{theorem}
\label{T0_1}
The functions $u=u_{\mathcal Re}(x,t)$ and $u_{\mathcal Im}(x,t)$ defined by 
\begin{eqnarray} 
\label{0.8K-G}
\hspace{-0.5cm} u_{\mathcal Re}(x,t)  
\!\!& \!\!= \!\!&\!\!
\int_{ t_0}^{t} db 
  \int_{ 0}^{t-b }   J_0\left(m\sqrt{(t-b)^2-r^2}\right)   w(x,r;b )  dr  
, \,\, x \in {\mathbb R}^n,  \\
\hspace{-0.5cm} u_{\mathcal Im}(x,t)  
\!\!& \!\!= \!\!&\!\!
\int_{ t_0}^{t} db 
  \int_{ 0}^{t-b }   I_0\left(m\sqrt{(t-b)^2-r^2} \right)  w(x,r;b )  dr  
, \,\, x \in {\mathbb R}^n,   
\end{eqnarray}
$m =|M|$, are solutions of the problems 
\begin{eqnarray*}  
u_{tt}-\Delta u +M^2 u =f(x,t) \,\,\mbox{\rm in}  \,\, {\mathbb R}^{n+1}, \quad
u(x,t_0)=0,\,\,\,  u_t(x,t_0)=0   \,\,\, \mbox{\rm in}  \,\, {\mathbb R}^{n}\,,  
\end{eqnarray*}
with $ M^2>0$ and $M^2<0 $, respectively. 
Here $w(x,t;b ) $ is a solution of  (\ref{wn}). 
\end{theorem}

\begin{definition} 
The integral operator  (\ref{Oper_K}) is said  to be a generator of the solution operator for some equation 
  if the operator $ G= {\mathcal K}\circ  {\mathcal WE} $   gives a solution operator for that equation.   
\end{definition}
\medskip

\noindent
{$2^0$ Tricomi-type equations.} The first example linkings  to the operator with the variable coefficient is generated by the kernel  $K (t;r,b)= E(0,t;r,b)$, where
the function $  E(x,t;r,b)$ \cite{YagTricomi} is defined by 
\begin{eqnarray} 
\label{E}
\hspace*{-0.5cm} E(x,t;r,b)
& := &
c_k \left(  (\phi (t)  + \phi (b))^2  -(x-r)^2 \right)^{-\gamma } \nonumber \\
&  &
\times  F \left(\gamma , \gamma ;1; \frac{(\phi (t)  - \phi (b))^2 - (x-r)^2}
{(\phi (t)  + \phi (b))^2 - (x-r)^2}  \right) \,,
\end{eqnarray}
with $\gamma := \frac{k}{2k+2} $, $c_k = (k+1)^{-k/(k+1)}2^{-1/(k+1)}$, $2k=l \in {\mathbb N}$, $x^2:=|x|^2$, and the distance function $\phi  =\phi (t)  $ is 
\begin{eqnarray}
\label{phi_Tric}
\phi (t)= \frac{1}{k+1}t^{k+1},\qquad 
\end{eqnarray}
while $F\big(a, b;c; \zeta \big) $ is the Gauss's hypergeometric function.  Here we assume that $2k \in{\mathbb N}$
but later on we consider the case of $l \in  {\mathbb R}$. It is proved in \cite{YagTricomi} that for an integer non-negative $l$,  for the smooth function $f=f(x,t)$, the function
\begin{eqnarray*}
u(x,t) 
& = &
2    c_l \int_{ 0}^{t} db 
  \int_{ 0}^{ \phi (t)- \phi (b)}  \left(  (\phi (t)  + \phi (b))^2  -r^2 \right)^{-\gamma } \\
&  &
\times F \left(\gamma , \gamma ;1; \frac{(\phi (t)  - \phi (b))^2 - r^2}
{(\phi (t)  + \phi (b))^2 - r^2}  \right)  w(x,r;b )  dr  
, \quad  t>0,
\end{eqnarray*}
 solves the Tricomi-type equation
\begin{eqnarray}
\label{Tric_eq} 
 u_{tt}-t^{l}\Delta u =f(x,t) \quad \mbox{\rm in}  \quad  {\mathbb R}_+^{n+1}:=\{(x,t)\,|\,  x \in {\mathbb R}^{n},\, t>0 \},   
\end{eqnarray}
and takes vanishing initial values
\begin{eqnarray}
\label{vanish_initial}  
u(x,0)=0,\quad  u_t(x,0)=0   \quad \mbox{\rm in}  \,\, {\mathbb R}^{n}\,. 
\end{eqnarray}

\noindent
{$3^0$ The wave equation in the Robertson-Walker spacetime: De~Sitter  spacetime.} The next interesting example we obtain if we set $K (t;r,b)= E(0,t;r,b)$, where
the function $  E(x,t;r,b)$ \cite{Yag_Galst_CMP} is defined by 
\begin{eqnarray} 
\label{E0_16}
\hspace*{-0.5cm} E(x,t;r,b)
& := &
\left(  (e^{-b}  + e^{-t})^2  -(x-r)^2 \right)^{-\frac{1}{2}}\nonumber \\
&  &
\times F \left(\frac{1}{2}, \frac{1}{2};1; \frac{(e^{-t}  - e^{-b})^2 - (x-r)^2}
{(e^{-t}  + e^{-b})^2 - (x-r)^2}  \right) \,,
\end{eqnarray}
and 
$
\phi (t):= e^{-t}  
$. 
It is proved in \cite{Yag_Galst_CMP} that  defined by the integral transform (\ref{Oper_K}) with the 
kernel (\ref{E0_16})  the function
\begin{eqnarray*} 
u(x,t)
& = &
 2\int_0^t \,db \int_{ 0}^{  e^{-b}-e^{-t}  } \left(  (e^{-b}  + e^{-t})^2  -r^2 \right)^{-\frac{1}{2}} \\
&  &
\times F \left(\frac{1}{2}, \frac{1}{2};1; \frac{(e^{-t}  - e^{-b})^2 - r^2}
{(e^{-t}  + e^{-b})^2 - r^2}  \right)w(x,r;\tau )\,dr  
\end{eqnarray*}
 solves the wave equation in the Robertson-Walker spaces arising in the de~Sitter 
model of the universe (see, e.g. \cite{Moller}),
\[ 
 u_{tt}-e^{-2t}\Delta u =f(x,t) \quad \mbox{\rm in}  \,\, {\mathbb R}_+^{n+1},  
\]
and takes vanishing initial data (\ref{vanish_initial}).

\noindent
{$4^0$  The wave equation in the Robertson-Walker spacetime: anti-de~Sitter  spacetime}. The third  example we obtain if we set $K (t;r,b)= E(0,t;r,b)$, where
the function $  E(x,t;r,b)$ is defined by (see \cite{Yag_Galst_JMAA})
\begin{eqnarray} 
\label{E_017}
\hspace*{-0.5cm} E(x,t;r,b)
& := &
\left(  (e^{b}  + e^t)^2  -(x-r)^2 \right)^{-\frac{1}{2}} \nonumber \\
&  &
\times F \left(\frac{1}{2}, \frac{1}{2};1; \frac{(e^t  - e^{b})^2 - (x-r)^2}
{(e^t  + e^{b})^2 - (x-r)^2}  \right) \,,
\end{eqnarray}
while the distance function is
$\phi (t):= e^t$ . 
In that case the function
$u=u(x,t)$  produced by the integral transform  (\ref{Oper_K}) with   $t_0=0$  and the 
kernel (\ref{E_017}),   solves the wave equation in the Robertson-Walker space arising in the anti-de~Sitter 
model of the universe (see, e.g. \cite{Moller}),
\[ 
 u_{tt}-e^{2t}\Delta u =f(x,t) \quad \mbox{\rm in}  \,\, {\mathbb R}_+^{n+1}.  
\]
Moreover, it takes vanishing initial values (\ref{vanish_initial}).
\medskip

\noindent
{$5^0$  The wave equation in the Einstein-de Sitter spacetime}. 
If we allow negative $l \in{\mathbb R}$ in (\ref{E}) and, in that way, simplify the Gauss's hypergeometric function of the kernel
of the integral transform, then we obtain another way to get new operators of the above described hierarchy.
In fact,  in the hierarchy of the  hypergeometric functions $ F \left(a, b;c; \zeta  \right)$ there  are functions which are polynomials. This is a case, in particular,
 of the parameter $a=-m$,  
where 
$m \in {\mathbb N}$. More precisely, if $m \in {\mathbb N}$, then $l=-4m/(2m+1) >-2$ and
$
  F \left(-m, -m;1; \zeta  \right) =  \sum_{n=0}^{m}\left( \frac{m(m-1)\cdots (m+1-n)}{n!}\right)^2 z^n 
$. 
In that case we choose the distance function $\phi (t)=(2m+1)t^{\frac{1}{2m+1}} $ and the kernel $K(t;r,b)$ as follows
\begin{eqnarray}  
K(t;r,b)  
&  = &
c_m  \sum_{n=0}^{m}\left( \frac{m(m-1)\cdots (m+1-n)}{n!}\right)^2 \nonumber \\
&  &  
\label{negat_l}
\times \left(  (2m+1)^2(t^{\frac{1}{2m+1}}   +  b^{\frac{1}{2m+1}} )^2  -r^2 \right)^{-m -n} \nonumber \\
&  &
\times \left(  (2m+1)^2(t^{\frac{1}{2m+1}}   -  b^{\frac{1}{2m+1}} )^2 - r^2 \right) ^n .
\end{eqnarray}
Thus the integral transform $\mathcal K$ allows us to write the representation for the  solution  of the equation
 \begin{eqnarray*}
\  u_{tt} 
-  t^{-\frac{4m}{2m+1}} \Delta   u  =f\quad \mbox{\rm in}  \quad  {\mathbb R}_+^{n+1}  \,.
\end{eqnarray*}
Moreover, in the hierarchy of the  hypergeometric functions
the  simplest non-constant function is  
$
  F \left(-1, -1;1; \zeta  \right) = 1+ \zeta $. 
The  exponent \,$   l$\, leading to   the  function $F \left(-1, -1;1; \zeta  \right) $   is exactly the exponent\, $   l=-4/3$ \, of the 
{\sl  wave equation }(and of the {\sl  metric tensor}) {\sl  in the Einstein~\&~de~Sitter spacetime}. In that case of $m=1$ the kernel  $K(t;r,b)$ of (\ref{negat_l}) is 
$
K(t;r,b)= 
\frac{1}{18} \left(   9t^{2/3}   + 9b^{2/3}    -r^2 \right)$. 
Consequently,  the function 
\begin{equation} 
\label{sol_EdS}
u(x,t)  
=  
   \int_{ 0}^{t} db 
  \int_{ 0}^{ 3t^{1/3}- 3b^{1/3}}   \frac{1}{18} \left(   (3t^{1/3})^2   + (3b^{1/3})^2     -r^2 \right)  w(x,r;b )  dr  
, 
\end{equation}
$x \in {\mathbb R}^n$, \,$t>0$, \, solves (see \cite{Galstian_Kinoshita_Yagdjian}) the equation
 \begin{eqnarray}
 \label{0.8}
u _{tt}
-  t^{-{4}/{3}}  \Delta   u  =f\quad \mbox{\rm in} \quad  {\mathbb R}_+^{n+1}  \,,
\end{eqnarray}
and  takes vanishing  initial data  (\ref{vanish_initial}) provided that the function $w $ is the image of $f$, that is $w ={\mathcal WE} (f) $. 
Because of the singularity in the coefficient of equation (\ref{0.8}), 
the Cauchy problem is not well-posed.  In order to obtain a  well-posed problem the initial conditions must be modified to the weighted initial value conditions.

In fact, the operator of equation (\ref{0.8}) coincides  with the principal part of the   
 wave equation  in the Einstein~\&~de~Sitter spacetime. 
We remind that the   Einstein\,\&\,de~Sitter model  (EdeS model) of the universe was first proposed jointly by Einstein and de~Sitter    in 1932. It is the simplest non-empty expanding model with
the line-element $ ds^2 = - dt^2 + a_0^2t^{4/3} \left( dx ^2+ dy ^2+dz^2 \right)$. 
The  covariant linear wave equation with the source term $f$ written in these coordinates  is
\[
\left(  \frac{\partial}{\partial t} \right)^2  \psi 
-   t^{-{4}/{3}}  \sum_{i=1,2,3}\left( \frac{\partial }{\partial x^i} \right)^2   \psi +   \frac{2 }{ t }     \frac{\partial}{\partial t}    \psi =f\,.
\]
The last equation belongs to the family of 
the non-Fuchsian  partial differential equations.  There is very advanced theory of such equations (see, e.g.,
\cite{Mandai,Tahara}), but according to  our knowledge the weighted initial value problem suggested in \cite{Galstian_Kinoshita_Yagdjian} (see (\ref{weighted_prob}) below) is the original one.
Assume that 
$ f(x,t) \in C^\infty ({\mathbb  R}^n\times (0,\infty))$, and that with some $\varepsilon >0$ one has
\[
|\partial_x^\alpha f(x,t)| + |t \partial_t \partial_x^\alpha f (x,t)| \leq C_\alpha t^{\varepsilon -2 }  \quad \mbox{for all}  \,\, x \in {\mathbb R}^n, 
\quad    \mbox{and for  small}\,\, t>0,
\]
and for every $\alpha  $, $|\alpha | \leq  [(n+1)/2]$. 
It is proved in \cite{Galstian_Kinoshita_Yagdjian}    that 
the function
\begin{eqnarray*}
\psi (x,t)
& = & 
\frac{1}{18t}  \int_{ 0}^{t} db
\int_{ 0 }^{ 3t^{1/3}-3b^{1/3} }   b  
 \big(9t^{2/3}  + 9b^{2/3}-r^2  \big) w (x,r ; b)\, dr 
\end{eqnarray*}
solves the problem 
\begin{eqnarray}
\label{weighted_prob}
\cases{ 
 \psi_{tt} - t^{-4/3}\bigtriangleup    \psi +   2   t^{-1}          \psi_t = f (x,t),  \qquad t>0 ,\,\, x \in {\mathbb R}^n,\cr 
 \displaystyle   \lim_{t\rightarrow 0}\, t \psi  (x,t) = 0, \quad   
\displaystyle   
\lim_{t\rightarrow 0} 
  \left(  t \psi_t  (x,t) + \psi  (x,t)  \right)
=  0,   \,\, x \in {\mathbb  R}^n , }  
\end{eqnarray}
provided that the function $w $ is the image of $f$, that is $w ={\mathcal WE }(f) $.

\section{The fundamental solutions of the operators} 
\setcounter{equation}{0}

In this section we apply the method from the previous section to construct the fundamental solutions of the operators  
for the above listed equations. They  are hyperbolic equations and therefore they have the fundamental solutions with the support in the forward or backward light cones.
First we consider the string equation. 
In order to write the  fundamental solution with the support in the forward  light cone we look for ${\mathcal E}\in {\mathcal D}'  ({\mathbb R}^{n+1})$  such that
\[ 
 {\mathcal E}_{tt}- {\mathcal E}_{xx}=\delta (x-x_0)\delta (t-t_0)  \quad \mbox{\rm in}  \,\, {\mathbb R}^{2},  
\quad \mbox{\rm supp}\, {\mathcal E}\subseteq \{(x,t)\,|\, t \geq t_0,\,\, x \in {\mathbb R}\}
\]
Then ${\mathcal E}(x,t;x_0,t_0) = {\mathcal E}(x-x_0,t-t_0;0,0) $ and if we denote by $ D(x_0,t_0)$ the forward  light cone 
$D(x_0,t_0):= \{(x,t)\in {\mathbb R}^{n+1}\,|\, |x-x_0| \leq (t-t_0)   \}$. 
It is well known that ${\mathcal E}(x,t;x_0,t_0) =  $ if $(x,t) \in D (x_0;t_0) $ and ${\mathcal E}(x,t;x_0,t_0) = 0$ otherwise.
 We follows approach of the previous section and rewrite this  fundamental solution   in the following way
\begin{eqnarray*} 
{\mathcal E}(x,t;x_0,t_0)  
& = &
H(t-t_0) \int_{0}^{t- t_0 } {\mathcal E}^{string}(x-x_0,z) \,dz \,,
\end{eqnarray*}
where $H(t-t_0)$ is the Heaviside  step function. The distribution 
\[
{\mathcal E}^{string} (x,t)= \frac{1}{2} \left\{ \delta (x+t) 
+ \delta (x- t) \right\}
\]
 is the fundamental solution of the Cauchy problem for the string equation:
\[
   {\mathcal E}_{tt}^{string}   -      {\mathcal E}_{xx}^{string} =0, \qquad  {\mathcal E}^{string}(x,0 )= \delta (x)
, \,\,\, {\mathcal E}^{string}_t  (x,0 )=0\,.
\] 
The string equation is partially, in direction of time, hypoelliptic, that implies
$
{\mathcal E}^{string}    \in C^\infty ( {\mathbb R}_t;{\mathcal D}'  ({\mathbb R}_x^{n})) 
$. 
Hence, for every test function  $\varphi \in C^\infty ({\mathbb R}^{n}) $, we have 
\begin{eqnarray*} 
<{\mathcal E}(x,t;\cdot ,t_0), \varphi (\cdot )>  
& = &
H(t-t_0)\int_{0}^{t- t_0 }< {\mathcal E}^{string}(x-\cdot ,z),\varphi (\cdot )> \,dz \,.
\end{eqnarray*}
Thus, the forward fundamental solution of the operator is given by the integral transform   of   the fundamental solution of the Cauchy problem for the wave equation
corresponding to  the first datum.

We can generate a class of equations which allows explicit representation formulas for the fundamental solutions.
Indeed, if  we provide the integral transform with a kernel as follows:
\[ 
{\mathcal E}(x,t;x_0,t_0)  
= 
H(t-t_0)\int_{ 0}^{ |\phi (t)- \phi (t_0)|}   K (t;r,b)  {\mathcal E}^{wave}(x-x_0,r)  dr  
\,,
\]
$x \in {\mathbb R}^n$, \, $t \in  {\mathbb R}$, \,then we get the representations for the  fundamental solutions of the wide class of partial differential equations equations. 
In particular, if we plug in the integral transform the kernels used in the previous examples, then we obtain the corresponding fundamental solutions with the support in the forward light cone.

\noindent
$1^0 $ Klein-Gordon equation in  the Minkowski spacetime. 
If we choose the kernel $K (t;r,b)$   (\ref{0.6})  and  choose the distance function as $\phi (t)=t $, 
then it can be easily verified (see Theorem~\ref{T0_FS} below) that the distribution   
\[ 
{\mathcal E}(x,t;x_0,t_0)  
=  
H(t-t_0)\int_{ 0}^{ t-b}   J_0\left(\sqrt{(t-b)^2-r^2}\right)  {\mathcal E}^{wave}(x-x_0,r)  dr  
\,, 
\]
$x \in {\mathbb R}^n$,  $t \in  {\mathbb R}$, \,is the forward fundamental solutions  for the Klein-Gordon operator with a positive mass equals to $ 1$ in the Minkowski spacetime,
\[  
\left(  \partial^2_{t}-\Delta  + 1\right) {\mathcal E}(x,t;x_0,t_0)  =\delta (x-x_0)\delta (t-t_0)  \,\, \mbox{\rm in}  \,\, {\mathbb R}^{n+1},  \quad
 \mbox{\rm supp}\, {\mathcal E}\subseteq D(x_0,t_0).  
\]
provided that $ {\mathcal E}^{wave}(x,t)$ is  the fundamental solution of the Cauchy problem corresponding to  the first datum with the support   at the origin, for 
the wave equation in the Minkowski spacetime. 
We emphasis that the distributions  ${\mathcal E}(x,t;x_0,t_0)$ and  $ {\mathcal E}^{wave}(x,t)$  solve {\sl different} equations.

If we now choose the kernel $K (t;r,b)$ (\ref{BesselI_0})
 and the distance function as $\phi (t)=t $, 
then  the distribution   
\[ 
{\mathcal E}(x,t;x_0,t_0)  
 =  
H(t-t_0)\int_{ 0}^{ t-t_0}   I_0\left(\sqrt{(t-b)^2-r^2}\right)  {\mathcal E}^{wave}(x-x_0,r)  dr  
\,, 
\] 
$x \in {\mathbb R}^n$,  $t \in  {\mathbb R}$, \,is the forward fundamental solutions  for the Klein-Gordon operator with an imaginary mass   in the Minkowski spacetime,
\[  
\left(  \partial^2_{t}-\Delta  - 1\right) {\mathcal E}(x,t;x_0,t_0)  =\delta (x-x_0)\delta (t-t_0)   \,\, \mbox{\rm in}  \, {\mathbb R}^{n+1},   
 \mbox{\rm supp}\, {\mathcal E}\subseteq D(x_0,t_0) .  
\]
The following theorem  can be easily proved by direct substitution.
\begin{theorem}
\label{T0_FS}
The distributions ${\mathcal E}_{\mathcal Re}(x,t;x_0,t_0)$ and ${\mathcal E}_{\mathcal Im}(x,t;x_0,t_0)  $ defined by 
\begin{eqnarray*} 
{\mathcal E}_{\mathcal Re}(x,t;x_0,t_0)  
\!\! &\!\!  = \!\! &\!\! 
H(t-t_0)\int_{ 0}^{ t-t_0}   J_0\left(m\sqrt{(t-b)^2-r^2}\right)  {\mathcal E}^{wave}(x-x_0,r)  dr  
,  \\
{\mathcal E}_{\mathcal Im}(x,t;x_0,t_0)  
\!\! &\!\!  = \!\! &\!\! 
H(t-t_0)\int_{ 0}^{ t-t_0}   I_0\left(m\sqrt{(t-b)^2-r^2}\right)  {\mathcal E}^{wave}(x-x_0,r)  dr ,    
\end{eqnarray*}
$x \in {\mathbb R}^n$,   $t \in  {\mathbb R}$,   are forward fundamental solutions  for the Klein-Gordon operators with a real and an imaginary mass  
\begin{eqnarray*} 
 \partial^2_{t}-\Delta  +M^2  \quad \mbox{\rm in}  \,\, {\mathbb R}^{n+1}, 
\end{eqnarray*}
with $ M^2 >0$ and $M^2<0 $, respectively. 
Here $m=|M|\geq 0$ and $ {\mathcal E}^{wave}(x,t)$ is  the fundamental solution of the Cauchy problem corresponding to  the first datum with the support   at the origin, for 
the wave equation  in the Minkowski spacetime.  
\end{theorem}

\noindent
$2^0 $ Tricomi-type equations. If we now choose the kernel $K (t;r,b)$ (\ref{E})
 and the distance function as (\ref{phi_Tric}), 
then it is proved in \cite{YagTricomi} that   
  ${\mathcal E}(x,t;x_0,t_0)$ is the  forward fundamental solution for the Tricomi-type equation (\ref{Tric_eq}):   
\begin{eqnarray*}  
{\mathcal E}(x,t;x_0,t_0)   
& = &
2    c_l H(t-t_0)
  \int_{ 0}^{ \phi (t)- \phi (t_0)}  \left(  (\phi (t)  + \phi (t_0))^2  -r^2 \right)^{-\gamma } \\
&  &
\times F \left(\gamma , \gamma ;1; \frac{(\phi (t)  - \phi (t_0))^2 - r^2}
{(\phi (t)  + \phi (t_0))^2 - r^2}  \right)  {\mathcal E}^{wave}(x-x_0,r)   dr,   
\end{eqnarray*}
$x \in {\mathbb R}^n$, \,   $ t_0 \geq 0$,\, with the support in the forward light cone
\[
D(x_0,t_0):= \{(x,t)\in {\mathbb R}^{n+1}\,|\, |x-x_0| \leq  \phi (t)- \phi (t_0)    \}\,.
\]

\noindent
$3^0 $ 
Klein-Gordon equations in the Robertson-Walker spacetime. The integral transform and, in particular, its kernel and the Gauss's hypergeometric function, open a way to establish a bridge between the wave equation (massless equation) and the Klein-Gordon equation (massive equation) in the curved spacetime.
 Indeed, if we allow the parameter $\gamma  $ 
of the function $F  (\gamma , \gamma ;1; z  ) $ to be a complex number, $\gamma  \in {\mathbb C}$, then this continuation into the complex plane produces the fundamental solutions 
${\mathcal E}_+(x,t;x_0,t_0) $ for the 
 Klein-Gordon operator in the de Sitter spacetime 
 as follows 
\begin{eqnarray*}
\hspace{-0.8cm} &  &
{\mathcal E}_+(x,t;x_0,t_0) 
 =   
2H(t-t_0) \int_{0}^  { e^{-t_0}- e^{-t}}   (4e^{-t_0-t })^{iM} \left((e^{-t_0 }+e^{-t})^2 - r^2\right)^{-\frac{1}{2}-iM    }   \\
\hspace{-0.8cm}&  &\hspace{2.2cm}
\times F\Big(\frac{1}{2}+iM   ,\frac{1}{2}+iM  ;1; 
\frac{ ( e^{-t_0}-e^{-t })^2 -r^2 }{( e^{-t_0}+e^{-t })^2 -r^2 } \Big){\mathcal E}^{wave} (x-x_0,r )  \, dr , 
\end{eqnarray*}
where the distribution ${\mathcal E}^{wave} (x,t )  $ is the fundamental solution of the Cauchy problem for the wave equation, while  the non-negative 
{\it curved mass} \,$M \geq 0$\, is defined as follows:
$M^2:=  \frac{n^2}{4} - m^2\geq 0$. 
The parameter $m$ is  mass of  particle. The fundamental solution ${\mathcal E}_-(x,t;x_0,t_0) $ with the support in the backward  light cone  admits  a similar representation.
 The fundamental solutions  ${\mathcal E}_+(x,t;x_0,t_0) $ and  ${\mathcal E}_-(x,t;x_0,t_0) $ are constructed in \cite{Yag_Galst_CMP} for the case of the large masses\, $m  \geq n /2$. 
The integral makes sense in the topology of the space of distributions.
The fundamental solutions for the  Klein-Gordon operator in the anti-de Sitter spacetime 
can be obtained by time inversion, $t \rightarrow -t$, from the fundamental solutions for the  Klein-Gordon operator in the de Sitter spacetime.

Moreover, 
the analytic continuation of this distribution   in parameter $M $ into ${\mathbb C}$  allows us to use it 
 also in the case of small mass\, $0 \leq m  \leq n /2$. Corresponding equation
\[
u_{tt} - e^{-2t} \bigtriangleup u  - M^2 u=   0,
\]
can be regarded as Klein-Gordon  equation with an imaginary mass. Equations with imaginary mass appear in 
several physical  models such as \, $\phi ^4$ \, field model, tachion (super-light) fields,  Landau-Ginzburg-Higgs equation  and others.

 More precisely, for small mass\, $0 \leq m  \leq n /2$
we  define the distribution    \,${\mathcal E}_+(x,t;x_0,t_0)$\, \,by 
 \begin{eqnarray*}
 &  &
{\mathcal E}_+(x,t;x_0,t_0) \\
& = & 
2H(t-t_0) \int_{0}^  { e^{-t_0}- e^{-t}}   (4e^{-t_0-t })^{-M} \Big((e^{-t_0 }+e^{-t})^2 - r^2\Big)^{-\frac{1}{2}+M    }  \\
&  &
\hspace{0cm} \times F\Big(\frac{1}{2}-M   ,\frac{1}{2}-M  ;1; 
\frac{ ( e^{-t_0}-e^{-t })^2 -r^2 }{( e^{-t_0}+e^{-t })^2 -r^2 } \Big){\mathcal E}^{wave} (x-x_0,r )  \, dr . 
\end{eqnarray*}

\noindent
$4^0 $ The above listed examples hint at some necessary condition on the pair  $ \phi $, $K (t;r,b)$  in order for that pair to produce a generator of the solution operator for some partial differential
equation.

\begin{theorem}
Assume that the integral transform (\ref{Oper_K}) with the distance function $ \phi $ and the kernel $K (t;r,b)$ generate  fundamental solution 
\[ 
{\mathcal E}(x,t;x_0,t_0)  
 =  
H(t-t_0)\int_{ 0}^{  |\phi (t)- \phi (t_0)| }   K (t;r,b)  {\mathcal E}^{wave}(x-x_0,r)  dr  
\,,
\]
$x \in {\mathbb R}^n$, $t \in  {\mathbb R}$, of the partial differential
equation
\begin{eqnarray}
\label{PDE}
u_{tt}- \sum_{i,j=1}^{n} (a_{ij}(x,t) u_{x_i})_{x_j} + \sum_{i=1}^{n} (b_{i}(x,t) u) _{x_i} +b(t)u_t +c(t)u=f 
\end{eqnarray}
with the real-analytic coefficients. Denote $V_1=V_1(t)$ and $V_2=V_2(t)$ two linearly independent solutions of the ordinary differential equation
\[
V'{}'  +b(t)V' +c(t)V =0, \qquad V_1(0)=1=V_2'(0),\quad V_1'(0)=0=V_2(0)\,.
\]

Then the function $K (t;r,b)$ satisfies identity 
\begin{eqnarray} 
\label{}
2      \int_{ 0}^{  |\phi (t)- \phi (b)| }   K (t;r,b) dr  = \frac{V_1(b) V_2(t)- V_1(t) V_2(b)}{V_1(b) V_2'(b)- V_1'(b) V_2(b )}
\end{eqnarray}
for all \, $t>b>0$. 
\end{theorem}

\noindent{\bf Proof.}
For every  function $f \in $ $C^\infty ({\mathbb R}\times [0,\infty))$, which 
for any given instant $t \ge 0$ has a compact support in $x$, the function
\begin{eqnarray*}
\hspace*{-0.0cm} v(x,t)  & = &
 \int_{ 0}^{t} db \int_{ 0}^{  |\phi (t)- \phi (b) |}   K (t;r,b)  w(x,r;b)  dr  \,    \\
\hspace*{-0.3cm} &  & \hspace*{0.5cm} 
\nonumber 
\end{eqnarray*}
where   \,$ w(x,r;b)= <{\mathcal E}^{wave}(x-\cdot ,r),f(\cdot ,b)>$,\, solves the equation (\ref{PDE})
and takes vanishing initial data.  It follows
\begin{eqnarray*}
V(t)& :=&
\int_{{\mathbb R}^n}  v(x,t) dx 
= 
\int_{ 0}^{t} db \int_{ 0}^{  |\phi (t)- \phi (b)| }   K (t;r,b) \left( \int_{{\mathbb R}^n}   w(x,r;b)  dx \right)\, dr  \,.
\end{eqnarray*}
On the other hand,
\begin{eqnarray*}
 \int_{{\mathbb R}^n}   w(x,r;b) \, dx = F(b), \quad F(b):=\int_{{\mathbb R}^n} f(x,b) \, dx\,.
\end{eqnarray*}
Hence,
\begin{eqnarray*}
V(t) 
& = & 
\int_{ 0}^{t} F(b)  \left( \int_{ 0}^{ | \phi (t)- \phi (b)| }   K (t;r,b)\, dr \right)db  \,.
\end{eqnarray*}
At the same time from the equation (\ref{PDE}) we obtain
\begin{eqnarray*} 
 \frac{d^2}{dt^2}V(t) +b(t) \frac{d }{dt }V(t) +c(t)V(t)=F(t) \,.
\end{eqnarray*}
Hence, 
\[
V(t)= \int_0^t F(b) \frac{V_1(b) V_2(t)- V_1(t) V_2(b)}{V_1(b) V_2'(b)- V_1'(b) V_2(b )}\, db\,.
\]
Thus, for the arbitrary function $f \in C^\infty ({\mathbb R}\times [0,\infty))$ for all $t$ one has 
\[
\int_{ 0}^{t}\!\Big( \int_{{\mathbb R}^n} f(x,b)  dx \Big)  \!
\Big( \int_{ 0}^{ | \phi (t)- \phi (b)| }   K (t;r,b)  dr - \frac{V_1(b) V_2(t)- V_1(t) V_2(b)}{V_1(b) V_2'(b)- V_1'(b) V_2(b )}\Big)db =0. 
\]
The theorem is proven.
\hfill $\Box$ 

\begin{corollary} 
1) For the Tricomi-type equation and for the wave equations in the de~Sitter and anti-de~Sitter spacetime    the following identities hold 
\begin{eqnarray*} 
t-b 
\!\! & \!\!=  \!\!&\!\!
2    c_k 
  \int_{ 0}^{ \frac{t^{k+1}}{k+1}  - \frac{b^{k+1}}{k+1} }  \left(  \left( \frac{t^{k+1}}{k+1}  + \frac{b^{k+1}}{k+1}  \right)^2  -r^2 \right)^{- \frac{k}{2k+2}  }\\
\!\!&\!\!  \!\!&\!\!
  \times 
  F \left(\frac{k}{2k+2} , \frac{k}{2k+2} ;1; \frac{\left(\frac{t^{k+1}}{k+1}  - \frac{b^{k+1}}{k+1} \right)^2 - r^2}
{\left(\frac{t^{k+1}}{k+1}  + \frac{b^{k+1}}{k+1} \right)^2 - r^2}  \right)     dr   \,,  \quad \,\, t > b \geq 0, \\ 
  t-b
\!\!&\!\! =  \!\!&\!\!
 2      \int_{ 0}^{   e^{-b}-e^{-t} }\left(  (e^{-b}  + e^{-t})^2  -r^2 \right)^{-\frac{1}{2}} F \left(\frac{1}{2}, \frac{1}{2};1; \frac{(e^{-t}  - e^{-b})^2 - r^2}
{(e^{-t}  + e^{-b})^2 - r^2}  \right) dr 
,\\
t-b
\!\!&\!\! =  \!\!&\!\!
2      \int_{ 0}^{   e^{t}-e^{b}  }  \left(  (e^{b}  + e^{t})^2  -r^2 \right)^{-\frac{1}{2}} F \left(\frac{1}{2}, \frac{1}{2};1; \frac{(e^{t}  - e^{b})^2 - r^2}
{(e^{t}  + e^{b})^2 - r^2}  \right) dr 
,  
\end{eqnarray*}
$ t > b$, respectively.\\
2) For the Klein-Gordon equation in the de~Sitter spacetime
   we have $b(t)=0$, $c(t)=M^2$ and the following identities hold. If the mass term $M^2$ is positive, then
\begin{eqnarray*}
\hspace{-0.5cm}\frac{1}{M} \sin M(t-b) 
& = & 
2\int_{0}^  { e^{-b}- e^{-t}}   (4e^{-b-t })^{iM} \Big((e^{-b }+e^{-t})^2 - r^2\Big)^{-\frac{1}{2}-iM    }   \\
\hspace{-0.5cm} &  & \hspace{1cm} \times F\Big(\frac{1}{2}+iM   ,\frac{1}{2}+iM  ;1; 
\frac{ ( e^{-b}-e^{-t })^2 -r^2 }{( e^{-b}+e^{-t })^2 -r^2 } \Big)  dr 
\,.
\end{eqnarray*}
If the mass term $M^2$ is negative, then
\begin{eqnarray*}
\hspace{-0.3cm}\frac{1}{|M|} \sinh |M|(t-b )
& = & 
2\int_{0}^  { e^{-b}- e^{-t}}   (4e^{-b-t })^{-|M|} \Big((e^{-b }+e^{-t})^2 - r^2\Big)^{-\frac{1}{2}+|M|    }   \\
\hspace{-0.3cm} &  & \hspace{0.5cm} \times F\Big(\frac{1}{2}-|M|   ,\frac{1}{2}-|M|  ;1; 
\frac{ ( e^{-b}-e^{-t })^2 -r^2 }{( e^{-b}+e^{-t })^2 -r^2 } \Big)   dr  
\,. 
\end{eqnarray*}
\end{corollary}

These identities   were used in \cite{YagTricomi_GE} and \cite{yagdjian_DCDS}  to prove a blow up phenomenon for the semilinear Tricomi-type  equation and the Klein-Gordon equation
in the de~Sitter spacetime.

\section{Estimate of the tail inside of  the light cone. De Sitter spacetime}
\setcounter{equation}{0}

In this section we consider the equation
\begin{equation}
\label{NWE} 
u_{tt} - e^{-2t} \Delta  u=  0 \,. 
\end{equation}
The forward and backward fundamental solutions for the operator of the last equation is constructed in \cite{Yag_Galst_CMP}. By means of those fundamental solutions
 the fundamental solutions of the Cauchy problem are given as the Fourier integral operators 
  in the  domain of hyperbolicity, $t>0$, 
with the data prescribed at $t=0$,
\begin{equation}
\label{CD} 
u(x,0)=\varphi _0 (x), \quad u_t(x,0)=\varphi _1 (x), \quad x\in {\mathbb R}^n \,.
\end{equation}
 The formula for the solution of this problem is given by Theorem~0.6~\cite{Yag_Galst_CMP} with $M=0$.
More precisely, the solution   of  equation (\ref{NWE}) with the initial data \, $ \varphi_0 $,  $ \varphi_1 \in C_0^\infty ({\mathbb R}^n) $, 
prescribed at $t=0$ is:  
\begin{eqnarray}
\label{1.30}
\hspace{-0.2cm} u(x,t) 
& = &
 e ^{\frac{t}{2}} v_{\varphi_0}  (x, \phi (t))
+ \, 2\int_{ 0}^{1} v_{\varphi_0}  (x, \phi (t)s) K_0(\phi (t)s,t)\phi (t)\,  ds  \nonumber \\
& &
+\, 2\int_{0}^1   v_{\varphi _1 } (x, \phi (t) s) 
  K_1(\phi (t)s,t) \phi (t)\, ds 
, \quad x \in {\mathbb R}^n, \,\, t>0,  
\end{eqnarray}
where $\phi (t):= 1-e^{-t} $. 
The function $v_\varphi  (x, \phi (t) s)$  coincides with the value $v(x, \phi (t) s) $ 
of the solution $v(x,t)$ of the Cauchy problem
$v_{tt}-  \bigtriangleup v =0$, \, $ v(x,0)= \varphi (x)$, \, $ v_t(x,0)=0$. 
The kernels  $K_0(z,t)   $    and $K_1(z,t)   $ are defined by  
$ K_0(z,t) $ $
 := 
- \left[  \frac{\partial }{\partial b}   E(z,t;0,b) \right]_{b=0} $ and $K_1(z,t) 
 := 
  E(z ,t;0,0) $, respectively, 
 where $E(x,t;r,b) $ is given by (\ref{E0_16}). Thus,  
\begin{eqnarray*}
\hspace{-0.5cm} &  &
 K_0(z,t) \\
\hspace{-0.5cm} &  = &   
\frac{1}{ [(1-e^{ -t} )^2 -  z^2]\sqrt{(1+e^{-t } )^2 - z^2} }\Bigg[  \big(  e^{-t} -1  \big) 
F \Big(\frac{1}{2}   ,\frac{1}{2}  ;1; \frac{ ( 1-e^{-t })^2 -z^2 }{( 1+e^{-t })^2 -z^2 }\Big) \\
\hspace{-0.5cm} &  &
  +   \big( 1-e^{-2 t}+  z^2 \big) \frac{1}{2}
F \Big(-\frac{1}{2}   ,\frac{1}{2}  ;1; \frac{ ( 1-e^{-t })^2 -z^2 }{( 1+e^{-t })^2 -z^2 }\Big) \Bigg],\, 0\leq z <  1-e^{-t},
\end{eqnarray*} 
and
\[ 
K_1(z,t)  
  =  
\big((1+e^{-t })^2 -   z  ^2\big)^{-\frac{1}{2} } 
F\left(\frac{1}{2}  ,\frac{1}{2}  ;1; 
\frac{ ( 1-e^{-t })^2 -z^2 }{( 1+e^{-t })^2 -z^2 } \right)\!, \, 0\leq z\leq  1-e^{-t}. 
\]

It is important that the formula (\ref{1.30}) 
regarded in the topology of the continuous functions of variable $t$ with the values in the 
distributions space  ${\mathcal D}'( {\mathbb R}^n) $,
is applicable to the distributions $\varphi_0, \varphi_1 \in {\mathcal D}'( {\mathbb R}^n) $ 
as well.

 Recall that a wave equation is said to satisfy Huygens'
principle if the solution
vanishes at all points which cannot be reached from the initial data by a null geodesic, that is   
there is no tail. An exemplar equation satisfying Huygens' principle is the wave
equation in $n + 1$ dimensional Minkowski spacetime for odd $n \geq  3$. According
to Hadamard's conjecture  this is the only (modulo transformations of coordinates and unknown function) huygensian
linear second-order hyperbolic equation.  
Counterexamples to Hadamard's conjecture, which have been found, do not change the fact that Huygens' property is a
very rare and unstable, with respect to the perturbations, phenomenon. It is natural, therefore,  to ask if there are other hyperbolic second-order  equations which preserve Huygens' property approximately, in the sense that
the tail which is left behind the wave front is comparatively  small. For the equation in the de~Sitter spacetime Huygens' principle is not valid 
\cite{Yag_Galst_CMP}. 
In this section we show the way how our approach can be applied to derive the pointwise estimates of the tail.

We call the ``tail'' the part containing the integrals of (\ref{1.30}), that is
\begin{eqnarray*}  
T (x,t) 
&   :=  &
 2\int_{ 0}^{1} v_{\varphi_0}  (x, \phi (t)s) K_0(\phi (t)s,t)\phi (t)\,  ds  \nonumber \\
& &
+\, 2\int_{0}^1   v_{\varphi _1 } (x, \phi (t) s) 
  K_1(\phi (t)s,t) \phi (t)\, ds 
, \quad x \in {\mathbb R}^n, \,\, t>0  
\,.
\end{eqnarray*}
Hence,  
$
T (x,t)   =  
u(x,t) -    e ^{\frac{t}{2}} v_{\varphi_0}  (x, \phi (t))$. 
The tail is   of considerable interest in many aspects in the physics, and in particular, in the  General Relativity  \cite{Bizon},  \cite{Sonego-Faraoni}.

In this section we restrict ourselves to the case of one-dimensional $x$. For the one-dimensional wave equation in the Minkowski spacetime Huygens'
principle is not valid, and, consequently, one can not anticipate it for one-dimensional equation in the de~Sitter spacetime. But the last one   reveals all difficulties and  technical details 
allowing to overcome  those   difficulties in the  case of  $x\in {\mathbb R}^n$ with $n \geq 2$. The results for the case of $n \geq 2$ will be published in a forthcoming paper.
We start with simple example.

\noindent
{\bf Example.}
 Let  $\varphi _0 (x)   = H(x)$, $\varphi _1 (x)  = 0$,
where $H$ is the Heaviside step function. Then $u(x,t)   
   =   
\frac{1}{2} e ^{\frac{t}{2}}  \big[ 
H   (x+ 1-e^{-t})  
+     H   (x -  1 +e^{-t})  \big]  
+ \int_{ 0}^{1-  e^{-t}} \big[ 
H   (x - z)  
+     H   (x  + z)  \big] K_0(z,t)\,  dz$
and, consequently,
$
T(x,t)
 =  
\int_{ 0}^{1-  e^{-t}} 
H   (x - z)  K_0(z,t)\,  dz
+  \int_{ 0}^{1-  e^{-t}} H   (x  + z)K_0(z,t)\,  dz 
$.
 Consider the point  $x_0$ such that $0 \leq x_0 <1-e^{-t}$. Then  we have 
\[
u(x_0,t)   
   =   
\frac{1}{2} e ^{\frac{t}{2}}    
+ \int_{ 0}^{x_0}K_0(z,t)\,  dz  
+ \int_{ 0}^{1-  e^{-t}}K_0(z,t)\,  dz = \frac{1}{2} + \int_{ 0}^{x_0}K_0(z,t)\,  dz
\]
while $
T(x_0,t)=\frac{1}{2}- \frac{1}{2} e ^{\frac{t}{2}}+ \int_{ 0}^{x_0}K_0(z,t)\,  dz= u(x_0,t)  - \frac{1}{2} e ^{\frac{t}{2}} $. In particular,
\begin{eqnarray*}
&  &
\frac{|T(x_0,t)|}{|u(x_0,t)-  T(x_0,t)|} =1- e^{-t/2}-2 e^{-t/2} \int_{ 0}^{x_0}K_0(z,t)\,  dz  \leq 2(1- e^{-t/2} ),\\
&  &
\lim_{t\to \infty} \lim_{x \to (1- e^{-t})^-}\frac{|T(x,t)|}{|u(x,t)-  T(x,t)|} =2.
\end{eqnarray*}
Thus, the tail dominates  
the  huygensian part of the solution.  If  $  x_0 >1$, then $u(x_0,t)= 1$, while $T(x_0,t)= 1-  e^{t/2}$. 
 \smallskip

Suppose now that initial data $\varphi _0$ and $\varphi _1$ are the
homogeneous  functions,
\begin{equation}
\label{idprop} 
  \varphi _0 (x)   = C_0 |x|^{-a}, \quad  \varphi _1 (x)  = C_1 |x|^{-b} \,.
\end{equation} 
The next theorem gives the pointwise estimate for the tail.

\begin{theorem}  
\label{TLin} 
Consider the Cauchy problem for the equation (\ref{NWE}), (\ref{CD}) with $n=1$ and the   data
  (\ref{idprop}), where  $a,b \in (1/2,1)$. Then in the light cone emitted   by the origin, that is on the set $
\{(x,t)\,|\, |x| < 1- e^{-t},\,\, t \geq 0 \} $, the solution, and, consequently, the tail,  satisfy 
\begin{eqnarray*}
| T(x,t)|
& = &
\left| u(x,t) -  \frac{1}{2}C_0 e ^{\frac{t}{2}}  \Big[ |x+ 1-e^{-t}|^{-a}  
+     |x -  1 +e^{-t}|^{-a}  \Big]  \right| \\
& \leq  &
|C_0| C\left(1+ t     \right)e^{\frac{1}{2}t} e^{at}  (1+e^{t}(1-|x|))^{1/2-a} \\
&  &
+|C_1| C\left(1+ t     \right) e^{-\frac{1}{2}t} e^{bt}  (1+e^{t}(1-|x|))^{1/2-b}\,.
\end{eqnarray*}
\end{theorem}

\noindent{\bf Proof.}
The representation of the solution of the Cauchy problem for the one-dimensional case ($n=1$) 
of   equation (\ref{NWE})  is given by Theorem 0.4\cite{Yag_Galst_CMP} with $M=0$. 
More precisely,  the solution $u=u (x,t)$ of the Cauchy problem 
\[
u_{tt} - e^{-2t}u_{xx}  =0\, ,\qquad u(x,0) = \varphi_0  (x) \,, \qquad u_t(x,0) =\varphi_1  (x)\,  ,
\]
with \,$\varphi_0    ,  \varphi_1  $ of (\ref{idprop}),   can be represented as follows
\begin{eqnarray*}
u(x,t)   
&  =   &
\frac{1}{2} e ^{\frac{t}{2}}  \Big[ 
\varphi_0   (x+ 1-e^{-t})  
+     \varphi_0   (x -  1 +e^{-t})  \Big]\\
&  & 
+ \int_{ 0}^{1-  e^{-t}} \big[ 
\varphi_0   (x - z)  
+     \varphi_0   (x  + z)  \big] K_0(z,t)\,  dz \nonumber \\
&  &    
+ \,\,\int_{0}^{1-  e^{-t}} \,\Big[      \varphi_1    (x- z)  +   \varphi _1   (x + z)    \Big] K_1(z,t) dz  \nonumber
\,.
\end{eqnarray*}
Then,  it is evident that the solution   of the Cauchy problem with $C_1=0$ is 
\begin{eqnarray*}
 u(x,t)     
&  =    &
C_0 e ^{\frac{t}{2}}  \Big[ |x+ 1-e^{-t}|^{-a}  
+     |x -  1 +e^{-t}|^{-a}  \Big]  \\
&  & 
+ C_0\int_{ 0}^{1-  e^{-t}} \big[ 
|x - z|^{-a}  
+    |x  + z|^{-a}  \big]  K_0(z,t) \,  dz \,.
\end{eqnarray*} 
Consider for $x\geq 0$ the first term in the last integral. It can be estimated as follows: 
\begin{eqnarray*}
&  &
\left| \int_{ 0}^{1-  e^{-t}}  
|x - y|^{-a}  
 K_0(y,t) \,  dy \right| 
 \leq  
\int_{ 0}^{1-  e^{-t}}  
|x - y|^{-a}  
  |K_0(y,t)| \,  dy  \\
&  =    &
\int_{ 0}^{1-  e^{-t}}  
|x - y|^{-a}  
\frac{1}{ [(1-e^{ -t} )^2 -  y^2]\sqrt{(1+e^{-t } )^2 - y^2} } \\
&  &
 \times  \Bigg|\Bigg[  \big(  e^{-t} -1  \big) 
F \Big(\frac{1}{2}   ,\frac{1}{2}  ;1; \frac{ ( 1-e^{-t })^2 -y^2 }{( 1+e^{-t })^2 -y^2 }\Big)
\\
&  &
\hspace{2cm} +   \big( 1-e^{-2 t}+  y^2 \big) \frac{1}{2}
F \Big(-\frac{1}{2}   ,\frac{1}{2}  ;1; \frac{ ( 1-e^{-t })^2 -y^2 }{( 1+e^{-t })^2 -y^2 }\Big) \Bigg]\Bigg| \,  dy.
\end{eqnarray*} 
If we denote $z=e^{t}$ and make a change $y=e^{-t}r $ in the last integral, then 
\begin{eqnarray}
\label{int_K_0}
&  &
\int_{ 0}^{1-  e^{-t}}  
|x - y|^{-a}  
  |K_0(y,t)| \,  dy \nonumber  \\
&  =    &
z^{ a }\int_{ 0}^{z-1 }  
|zx - r|^{-a}  
\frac{1}{   [(z- 1)^2 -  r^2]\sqrt{(z+1 )^2 - r^2} } \nonumber \\
&  &
   \Bigg|\Bigg[  \big(  z -z^2  \big) 
F \Big(\frac{1}{2}   ,\frac{1}{2}  ;1; \frac{ ( z-1)^2 -r ^2 }{( z+1)^2 - r ^2 }\Big)\nonumber \\
&  &
  +   \big( z^2-1+  r^2 \big) \frac{1}{2}
F \Big(-\frac{1}{2}   ,\frac{1}{2}  ;1; \frac{ ( z-1)^2 -r ^2 }{( z+1)^2 -r ^2 }\Big) \Bigg]\Bigg| \,  dr\,. 
\end{eqnarray} 
Define two zones $Z_1(\varepsilon, z) $ and $ Z_2(\varepsilon, z )$, 
as follows
 \begin{eqnarray}
\label{9.10}
Z_1(\varepsilon, z) 
& := &
\left\{ (z,r) \,\Big|\, \frac{ (z-1)^2 -r^2   }{ (z+1)^2 -r^2 } \leq \varepsilon,\,\, 0 \leq r \leq z-1 \right\} ,\\ 
\label{9.11}
Z_2(\varepsilon, z) 
& := &
\left\{ (z,r) \,\Big|\, \varepsilon \leq  \frac{ (z-1)^2 -r^2   }{ (z+1)^2 -r^2 },\,\, 0 \leq r \leq z-1  \right\},
\end{eqnarray}
respectively. Then we split integral of (\ref{int_K_0}) into the sum $I_1+I_2$ where: 
\begin{eqnarray*}
I_k & := &
 z^{ a }\int_{ Z_k(\varepsilon, z) }  
|zx - r|^{-a}  
\frac{1}{   [(z- 1)^2 -  r^2]\sqrt{(z+1 )^2 - r^2} } \\
&  &
\hspace{1cm}  \times   \Bigg|\Bigg[  \big(  z -z^2  \big) 
F \Big(\frac{1}{2}   ,\frac{1}{2}  ;1; \frac{ ( z-1)^2 -r ^2 }{( z+1)^2 - r ^2 }\Big)
 \\
&  &
\hspace{1.3cm}  +   \big( z^2-1+  r^2 \big) \frac{1}{2}
F \Big(-\frac{1}{2}   ,\frac{1}{2}  ;1; \frac{ ( z-1)^2 -r ^2 }{( z+1)^2 -r ^2 }\Big) \Bigg]\Bigg| \,  dr,\,\, k=1,2.
\end{eqnarray*} 
First, we  restrict ourselves to the   first zone $Z_1(\varepsilon, z) $. 
We follow the arguments have been used in the proofs  of Lemma~7.4~\cite{Yag_Galst_CMP} and Proposition~10.2~\cite{Yag_Galst_CMP}.
In the first zone we have
\begin{eqnarray*}
F\Big(\frac{1}{2} ,\frac{1}{2} ;1; \frac{ (z-1)^2 -r^2   }{ (z+1)^2 -r^2 }   \Big) 
\!\! & \!\! = \!\! &\!\! 
 1 +  \frac{1}{4}  \frac{ (z-1)^2 -r^2   }{ (z+1)^2 -r^2 }   
+ O\left(\left( \frac{ (z-1)^2 -r^2   }{ (z+1)^2 -r^2 }\right)^2\right)\!\!,  \\
F\Big(-\frac{1}{2} ,\frac{1}{2} ;1; \frac{ (z-1)^2 -r^2   }{ (z+1)^2 -r^2 }   \Big) 
\!\! & \!\! = \!\! &\!\! 
 1 -   \frac{1}{4}  \frac{ (z-1)^2 -r^2   }{ (z+1)^2 -r^2 }   
+ O\left(\left( \frac{ (z-1)^2 -r^2   }{ (z+1)^2 -r^2 }\right)^2\right)  \!\! . 
\end{eqnarray*}
We use the last formulas to estimate the term containing the hypergeometric functions:
\begin{eqnarray*}
\hspace{-0.5cm} &   &
\Bigg|  \big( z -z^{ 2 }    \big) 
F \Big(\frac{1}{2}   ,\frac{1}{2}  ;1; \frac{ ( z-1)^2 -r^2 }{( z+1)^2 -r^2 }\Big) \\
&  &  
\hspace{2cm} +    \big(z^{2 }- 1+  r^2  \big)  \frac{1}{2} 
F \Big(-\frac{1}{2}   ,\frac{1}{2}  ;1; 
\frac{ ( z-1)^2 -r^2 }{( z+1)^2 -r^2 }\Big)  \Bigg|  \nonumber \\
\hspace{-0.5cm} & \leq  & 
\frac{1}{2}\big[(z-1)^2 -r^2 \big] 
+\frac{1}{8}\frac{ (z-1)^2 -r^2   }{ (z+1)^2 -r^2 }\left| -(r^2+z^2-1  )+2(-z^2+z  )\right|  \nonumber  \\
\hspace{-0.5cm} &  &
 + \Big( \left|  z -z^{ 2 }  \right| 
+\left|  z^{2 }- 1+  r^2  \right|\Big)O\left(\left( \frac{ (z-1)^2 -r^2   }{ (z+1)^2 -r^2 }\right)^2\right) .
\end{eqnarray*}
In that zone   the integral can be estimated by:
\begin{eqnarray*}
I_1 & \leq &
C  \int_{ 0}^{z-1}  |zx - y|^{-a}
\Bigg[ \frac{1}{  \sqrt{(z+1)^2 - y^2} }   
\Bigg\{ 1+ z^2\frac{1 }{(z +1)^2 -y^2 }   \Bigg\} 
 \Bigg]  \,  d y  \,.
\end{eqnarray*}
The last integral can be written as follows:
\begin{eqnarray} 
\label{second}
&  &
\int_{ 0}^{z-1}  |zx - y|^{-a}
\Bigg[ \frac{1}{   \sqrt{z }\sqrt{z+1 - y } }   
\Bigg\{ 1+ z \frac{1 }{ z+1 -y  }   \Bigg\} 
 \Bigg]   d y  \nonumber  \\
&  &=
 \int_{ 1}^{z}  |z(1-x)-s|^{-a}\Bigg\{
  \frac{1}{\sqrt{z }}  (1+s)^{-1/2} 
+     \sqrt{z }
(1 +s )^{-3/2} \Bigg\} d s  .
\end{eqnarray}
Consider the first term of the last sum.  
Inside of the light cone $|x| < 1-e^{-t} $, that is, $|zx| < z-1  $,  and one has
\begin{eqnarray*} 
&  &
\frac{1}{\sqrt{z }} \int_{ 1}^{z}  |z(1-x)-s|^{-a}
 (1+s)^{-1/2} 
  \,  d s \\
& = &
\frac{1}{\sqrt{z }} \int_{ 1}^{z(1-x)}  (z(1-x)-s)^{-a}
(1+s)^{-\frac{1}{2}}  
  \,  d s \\
&  &+\frac{1}{\sqrt{z }} \int_{ z(1-x)}^{z}  (s-z(1-x))^{-a}
 (1+s)^{-\frac{1}{2}} 
  \,  d s   \,.
\end{eqnarray*}
On the other hand
\begin{eqnarray*} 
\hspace{-1cm} &  &
\frac{1}{\sqrt{z }} \int_{ 1}^{z(1-x)}  (z(1-x)-s)^{-a}
(1+s)^{-\frac{1}{2}}  
  \,  d s 
 = 
\frac{1}{\sqrt{z }} (1+z(1-x))^{-a} \\
\hspace{-1cm} &  &
\hspace{1.3cm} \times \left(\frac{\sqrt{\pi }\Gamma (1-a)}{\Gamma \left(\frac{3}{2}-a\right)} (1+z(1-x))^{1/2}  
-2 \sqrt{2}  F\left(\frac{1}{2},a;\frac{3}{2};\frac{2}{1+z(1-x)}\right)\right),
\end{eqnarray*}
while
\begin{eqnarray*} 
&  &
\frac{1}{\sqrt{z }} \int_{ z(1-x)}^{z} (s-z(1-x))^{-a}
(1+s)^{-\frac{1}{2}} 
  \,  d s  \\
& = &
-\frac{1}{\sqrt{z }}\frac{e^{i a \pi } \sqrt{\pi }  \Gamma (1-a)}{\Gamma \left(\frac{3}{2}-a\right)}(1+z(1-x))^{\frac{1}{2}-a}\\
&  &
+\frac{2e^{i a \pi }}{\sqrt{z }} \sqrt{1+z} \left( 1+z(1-x) \right)^{-a} F\left(\frac{1}{2},a,\frac{3}{2},\frac{1+z}{1+z(1-x)}\right).
\end{eqnarray*}
Thus,
\begin{eqnarray} 
&  &
\frac{1}{\sqrt{z }} \int_{ 1}^{z}  |z(1-x)-s|^{-a}
 \frac{1}{   \sqrt{1+s} }  
  \,  d s  \nonumber \\
 & \leq &
 z^{-\frac{1}{2}}(1+z(1-x))^{-a} \left(-\sqrt{\pi }\left(e^{i a \pi }-1 \right)  (1+z(1-x))^{\frac{1}{2}}\frac{ \Gamma (1-a)}{\Gamma \left(\frac{3}{2}-a\right)} \right.  \nonumber \\
 &  &
  -2 \sqrt{2} F\left(\frac{1}{2},a,\frac{3}{2},\frac{2}{1+z(1-x)}\right)\Bigg) \nonumber \\
&  &
\label{int1z}
+2 z^{-\frac{1}{2}}(1+z)^{\frac{1}{2}} (1+z(1-x))^{-a} e^{i a \pi }F\left(\frac{1}{2},a,\frac{3}{2},\frac{1+z}{1+z(1-x)}\right) .
\end{eqnarray}
Here the arguments of the hypergeometric functions satisfy  the inequalities $0 < \frac{2}{1+z(1-x)}<1$    and $ 1< \frac{1+z}{1+z(1-x)}$
for all $z(1-x)> 1$, $ x>0$, $ z>1$. Since $\Re (\frac{3}{2}-\frac{1}{2}-a)>0$ we have
\[
\left|F\left(\frac{1}{2},a,\frac{3}{2},\frac{2}{1+z(1-x)}\right) \right| \leq const \quad \mbox{\rm for all }\quad z(1-x)> 1  \,.
\]
Then, to estimate the last  term of (\ref{int1z}) we use (4) of Sec.2.10 \cite{B-E}: for all $z(1-x)\geq 1$, $z>1$, $x>0$,
\begin{eqnarray*} 
\hspace{-0.7cm}&  &
F\left(\frac{1}{2},a;\frac{3}{2};\frac{1+z}{1+z(1-x)}\right)  =
A_1\left(\frac{1+z}{1+z(1-x)}\right)^{-\frac{1}{2}}
+A_2\left(\frac{1+z}{1+z(1-x)}\right)^{-1}\\
\hspace{-0.7cm}&  &
\hspace{2.7cm} \times \left(1-\frac{1+z}{1+z(1-x)}\right)^{1-a} F\left(1,\frac{1}{2};2-a;1-\frac{1+z(1-x)}{1+z} \right) ,
\end{eqnarray*} 
where
$
A_1
=\frac{\sqrt{\pi } \Gamma (1-a)}{2 \Gamma \left(\frac{3}{2}-a\right)}$,\, $ 
A_2
=\frac{\sqrt{\pi }\Gamma  (a-1 )}{2 \Gamma (a)}$.  
Hence,
\begin{eqnarray*} 
&  &
z^{-\frac{1}{2}}(1+z)^{\frac{1}{2}} (1+z(1-x))^{-a}   F\left(\frac{1}{2},a;\frac{3}{2};\frac{1+z}{1+z(1-x)}\right) \\ 
& = &
A_1z^{-\frac{1}{2}}(1+z)^{\frac{1}{2}} (1+z(1-x))^{-a}  \left(\frac{1+z}{1+z(1-x)}\right)^{-\frac{1}{2}}
\\
&  &
+A_2z^{-\frac{1}{2}}(1+z)^{\frac{1}{2}} (1+z(1-x))^{-a}   \left(\frac{1+z}{1+z(1-x)}\right)^{-1}\\
&  &
\hspace{0.5cm} \times \left(1-\frac{1+z}{1+z(1-x)}\right)^{1-a}F\left(1,\frac{1}{2};2-a;1-\frac{1+z(1-x)}{1+z} \right).
\end{eqnarray*}
 Since $a>1/2$ one has
\begin{eqnarray*} 
&  &
\left|  F\left(1,\frac{1}{2};2-a;1-\frac{1+z(1-x)}{1+z} \right)\right| \leq const \left(\frac{1+z}{1+z(1-x)} \right)^{a-1/2}  
\end{eqnarray*}
for all $ z(1-x)\geq 1$. It follows
\begin{eqnarray*}  
&  &
\left|z^{-\frac{1}{2}}(1+z)^{\frac{1}{2}} (1+z(1-x))^{-a}   F\left(\frac{1}{2},a;\frac{3}{2};\frac{1+z}{1+z(1-x)}\right)\right| \\ 
& \leq  &
C z^{-\frac{1}{2}}  (1+z(1-x))^{1/2-a}      \,.
\end{eqnarray*}
Finally
\begin{eqnarray*} 
\frac{1}{\sqrt{z }} \int_{ 1}^{z}  |z(1-x)-s|^{-a}
 \frac{1}{   \sqrt{1+s} }  
  \,  d s 
& \leq  &
C z^{-\frac{1}{2}}  (1+z(1-x))^{1/2-a}   \,.   
\end{eqnarray*}
From this inequality now we derive estimate for the second integral of (\ref{second}),
\[ 
\sqrt{z }\int_{1}^{z}  |z(1-x)-s|^{-a}
\frac{1}{   (1 +s )^{3/2}}\,  d s 
 \leq  
C z^{\frac{1}{2}}  (1+z(1-x))^{1/2-a}  \,.  
\]

In the second zone we have
\[
\varepsilon \leq  \frac{ (z-1)^2 -r^2   }{ (z+1)^2 -r^2 } \leq 1 \quad \mbox{\rm and}  \quad 
\frac{ 1  }{ (z-1)^2 -r^2 }  \leq  \frac{ 1   }{ \varepsilon[(z+1)^2 -r^2] }\,.
\]
and due to 15.3.10 of Ch.15\cite{A-S} we obtain
\[
\left| F\Big(\frac{1}{2},\frac{1}{2};1; \frac{ (z-1)^2 -r^2   }{ (z+1)^2 -r^2 }    \Big) \right|   
  \leq   
C \left(1+ \ln  z     \right) 
, \quad \mbox{\rm for all}  \quad (z,r) \in Z_2(\varepsilon, z)  . 
\]
This allows us to  estimate  the integral over the second zone:
\begin{eqnarray*}
I_2 
& \leq &
\int_{ (z,r) \in Z_2(\varepsilon, z) }  
|zx - r|^{-a}  
\frac{z^2}{   [(z- 1)^2 -  r^2]\sqrt{(z+1 )^2 - r^2} } \left(1+ \ln  z     \right)\,  dr\\
& \leq &
 \left(1+ \ln  z     \right)\int_{ (z,r) \in Z_2(\varepsilon, z) }  
|zx - r|^{-a}  
\frac{z^2}{  ( (z+1 )^2 - r^2)^{3/2}}\,  dr\\
& \leq &
\left(1+ \ln  z     \right)\int_{  0}^{z-1}  
|zx - r|^{-a}  
\frac{z^2}{  ( (z+1 )^2 - r^2)^{3/2}} \,  dr\\
& \leq &
\left(1+ \ln  z     \right)\sqrt{z}\int_{  1}^{z }  
|z(1-x) - s|^{-a}  
  ( 1+s)^{-3/2} \,  ds\,.
\end{eqnarray*} 
Then from the last inequality and (\ref{}) we obtain,
\begin{eqnarray*}
I_2 & \leq &
C \left(1+ \ln  z     \right)z^{\frac{1}{2}}  (1+z(1-x))^{1/2-a}\,.
\end{eqnarray*} 
Now consider the case of $C_0=0$.
Then
\begin{eqnarray*} 
T(x,t)   
&  =   &
C_1\int_{ 0}^{1-  e^{-t}} \big[ 
|x - z|^{-b}  
+    |x  + z|^{-b}  \big]  K_1(z,t) \,  dz  
\,,
\end{eqnarray*}
and
\begin{eqnarray*} 
|T(x,t)|   
& \leq   &
C_1\int_{ 0}^{1-  e^{-t}} \big[ 
|x - z|^{-b}  
+    |x  + z|^{-b}  \big]  \big((1+e^{-t })^2 -   z  ^2\big)^{-\frac{1}{2} } \\
&  &
\hspace{2cm}\times 
\left| F\left(\frac{1}{2}  ,\frac{1}{2}  ;1; 
\frac{ ( 1-e^{-t })^2 -z^2 }{( 1+e^{-t })^2 -z^2 } \right)\right| \,  dz  \,.
\end{eqnarray*}
Consider case of $x \geq 0$. Then make change of variable $z=re^{-t}$:
\[ 
|T(x,t)|   
 \leq  
C_1e^{bt}\int_{ 0}^{e^{t}- 1 }  
|e^{t}x - r |^{-b}    \big((e^{t}+1)^2 -   r^2\big)^{-\frac{1}{2} } 
\left| F\left(\frac{1}{2}  ,\frac{1}{2}  ;1; 
\frac{ ( e^{t}-1)^2 -r^2 }{( e^{t}+1)^2 - r  ^2 } \right)\right|   dr  .
\]
Denote $z= e^{t}$. Then
\[
|T(x,t)|   
 \leq   
Cz^{b}\int_{ 0}^{z- 1 }  
|zx - r |^{-b}    \big((z+1)^2 -   r^2\big)^{-\frac{1}{2} } 
\left| F\left(\frac{1}{2}  ,\frac{1}{2}  ;1; 
\frac{ ( z-1)^2 -r^2 }{( z+1)^2 - r  ^2 } \right)\right| \,  dr .
\]
Since $1/2<b<1$, then for all $z>1$ the following estimate 
\begin{eqnarray*}
&  &
\int_{  0}^{ z  - 1} 
|zx - r |^{-b}((z  + 1)^2  - r^2  )^{-\frac{1 }{2}}  
F\left(\frac{1}{2},\frac{1}{2};1; 
\frac{ (z  - 1 )^2 -  r^2   }
{ (z  + 1 )^2  - r^2} \right)   d r \nonumber \\
& \leq &
C(1+ \ln   z )z^{-1/2}  (1+z(1-x))^{1/2-b}  
\end{eqnarray*}
is fulfilled. 
To prove the last estimate we rewrite the argument of the hypergeometric function as follows $\frac{ (z  - 1 )^2 -  r^2   }
{ (z  + 1 )^2  - r^2} = 1- \frac{  4z  }{ (z  + 1)^2 -r ^2}$.
If $
r \ge \sqrt{(z  + 1)^2-8z}$,  
then $
\frac{  4z  } { (z  + 1)^2 -r ^2}  \ge \frac{1}{2}$    and  $0< 1- \frac{  4z  } { (z  + 1)^2 -r ^2} \le \frac{1}{2}$
for such $r$ and $z$  imply
\[
\left| F\left(\frac{1}{2},\frac{1}{2};1; 1-
\frac{  4z  }
{ (z  + 1)^2 -r ^2} \right) \right| \le C  \,.         
\]
Hence we have
\begin{eqnarray*}
&  &
\int_{  \sqrt{(z  + 1)^2-8z}}^{ z  - 1} 
|zx - r |^{-b}((z  + 1)^2  - r^2  )^{-\frac{1 }{2}}  
F\left(\frac{1}{2},\frac{1}{2};1; 
1- \frac{  4z  }{ (z  + 1)^2 -r ^2} \right)   d r \\
&  & \leq
C\int_{  \sqrt{(z  + 1)^2-8z}}^{ z  - 1} 
|zx - r |^{-b}((z  + 1)^2  - r^2  )^{-\frac{1 }{2}}  
   d r \\
&  & \leq
Cz^{-1/2}\int_{  0}^{ z  - 1} 
|zx - r |^{-b}(z  + 1   - r  )^{-\frac{1 }{2}}  
   d r \\
&  & \leq
C z^{-1/2}  (1+z(1-x))^{1/2-b} .
\end{eqnarray*}
If $
r \le \sqrt{(z  + 1)^2-8z}$   and $z\ge 6$, 
then $8 <8z \leq (z  + 1)^2 - r^2 \leq (z  + 1)^2 $,
implies
\[
\left| F\left(\frac{1}{2},\frac{1}{2};1; 1-
\frac{  4z  }
{ (z  + 1)^2 -r ^2} \right) \right| \le C \left| \ln \left( \frac{  4z  }
{ (z  + 1)^2 -r ^2} \right)  \right|  \le C   (1+ \ln   z )  \,.       
\]
Hence 
\begin{eqnarray*}
\hspace{-0.5cm} &  &
\int_{0 }^{\sqrt{(z  + 1)^2-8z  }} 
|zx - r |^{-b}((z  + 1)^2  - r^2  )^{-\frac{1 }{2}}  
F\left(\frac{1}{2},\frac{1}{2};1; 
1- \frac{  4z  }{ (z  + 1)^2 -r ^2} \right)   d r \\
\hspace{-0.5cm} &  &\leq
\int_{0 }^{\sqrt{(z  + 1)^2-8z  }} 
|zx - r |^{-b}((z  + 1)^2  - r^2  )^{-\frac{1 }{2}}  
C   (1+ \ln   z )  d r \\
\hspace{-0.5cm} &  &\leq
C   (1+ \ln   z ) z^{-1/2} \int_{0 }^{z-1} 
|zx - r |^{-b}(z  + 1   - r )^{-\frac{1 }{2}}  
 d r \\
\hspace{-0.5cm} &  &\leq
C   (1+ \ln   z )z^{-1/2}  (1+z(1-x))^{1/2-b}\,.
\end{eqnarray*}
Theorem is proven. \hfill $\Box$

\begin{small}
\nocite{*}
\bibliographystyle{rendiconti}

\end{small}
\nocite{*}
\bibliographystyle{rendiconti}
\bibliography{yagdjian}
\vspace{-0.6cm}

\end{document}